\documentclass[11pt]{article}
\usepackage{amsfonts}
\newcommand{\eq}{\begin{equation}}
\newcommand{\en}{\end{equation}}
\newcommand{\giv}{\,|\,}
\newcommand{\re}[1]{\mbox{(\ref{#1})}}

\def\qed{\mbox{\rule{0.5em}{0.5em}}}

\newtheorem{Theorem}{Theorem}[section]
\newtheorem{theorem}[Theorem]{Theorem}
\newtheorem{lemma}[Theorem]{Lemma}
\newtheorem{corollary}[Theorem]{Corollary}
\newtheorem{construction}[Theorem]{Construction}
\newtheorem{proposition}[Theorem]{Proposition}
\newtheorem{example}[Theorem]{Example}
\newtheorem{defn}[Theorem]{Definition}

\newtheorem{question}[Theorem]{Question}
\newtheorem{conjecture}[Theorem]{Conjecture}
\newtheorem{condition}[Theorem]{Condition}
\newtheorem{remark}{Remark}
\newtheorem{problem}{Problem}

\def\endpf{\hfill $\qed$ \vskip .25in}

\newfont{\msbm}{msbm10}
\newfont{\eusb}{eusb10}
\newfont{\eusm}{eusm10}
\newfont{\eurb}{eurb10}
\newfont{\eurm}{eurm10}
\newfont{\eufb}{eufb10}
\newfont{\eufm}{eufm10}

\newcommand {\reals} {\mbox{\msbm\symbol{'122}}}

\newcommand{\te}{\rightarrow}
\newcommand{\ed}{\mbox{$ \ \stackrel{d}{=}$ }}

\newcommand{\eps}{\varepsilon}

\newcommand{\Lev}{L\'evy}

\newcommand{\hf} {{\mbox{${\textstyle\frac{1}{2}}$}}}

\newcommand{\lb}[1]{\label{#1}}

\newenvironment{thm}[1]{\begin{theorem}\label{#1}}{\end{theorem}}

\newenvironment{lmm}[1]{\begin{lemma}\label{#1}}{\end{lemma}}
\newenvironment{qst}[1]{\begin{question}\label{#1}}{\end{question}}
\newenvironment{cnj}[1]{\begin{conjecture}\label{#1}}{\end{conjecture}}

\newenvironment{crl}[1]{\begin{corollary}\protect\label{#1}}{\end{corollary}}

\newenvironment{exm}[1]{\begin{example}\protect\label{#1}}{\end{example}}
\newenvironment{prp}[1]{\begin{proposition}\protect\label{#1}}{\end{proposition}}

\def\omxxx{{\sf prob} (S^{d-1})}
\def\omxxxx{\mbox{{\sf cb-sets}}}
\def\Rd{{\mathbb R}^d}
\def\R{{\mathbb R}}
\def\F{{\cal F}}
\def\diseq{\ed}
\def\P{{\mathbb P}}
\def\E{{\mathbb E}}
\def\ball{{\cal B}}
\def\grid{{\sf GRID}}
\def\balls{{\sf BALLS}}
\def\cen{{\sf Cen}}

\def\range{{\sf range}}
\def\Xb{{\overline{X}}}
\def\Yb{{\overline{Y}}}
\def\|{\, | \,}

\def\ee{\epsilon}

\def\one{{\bf 1}}

\def\C{{\cal C}}
\def\edge{{\sf parent}}
\def\Set{{\sf set}}
\def\sph{\pi}
\def\cut{\, {\rm cut}}
\def\INT{\, {\sf interior}}
\def\rad{{\rm r}}

\begin{document}


\title{Where Did The Brownian Particle Go?}

\author{Robin Pemantle\thanks{Research supported in part by a Presidential
   Faculty Fellowship}, \, \,
Yuval Peres\thanks{Research supported in part by NSF grant DMS-9803597},
  \, \,
Jim Pitman\thanks{Research supported in part by NSF grant DMS-97-03961}, \, \,
Marc Yor \thanks{Research supported in part by NSF grant DMS-97-03961} \, \,
\\
\\
Technical Report No. 564
\\
\\
Department of Statistics\\
University of California\\
367 Evans Hall \# 3860\\
Berkeley, CA 94720-3860}
\date{December 3, 2000}
\maketitle
\begin{abstract}

Consider the radial projection onto the unit sphere
of the path a  $d$-dimensional Brownian motion $W$,
started at the center of the sphere and run for unit time. 
Given the occupation measure $\mu$ of this projected
path, what can be said about the terminal point $W(1)$, or about the
range of the original path?  In any dimension, for each Borel set $A \subseteq
S^{d-1}$, the conditional probability that the projection of 
$W(1)$ is in $A$ given $\mu(A)$ is just $\mu (A)$.  Nevertheless, 
in dimension $d \ge 3$, both the range and the terminal 
point of $W$ can be recovered with probability 1 from $\mu$. 
In particular, for $d \ge 3$ the conditional law of the projection of 
$W(1)$ given $\mu$ is not $\mu$.
In dimension~2 we conjecture that
the projection of $W(1)$ cannot be recovered almost surely from $\mu$,
and show that the conditional law of the projection of 
$W(1)$ given $\mu$ is not $\mu$.
\end{abstract}
\newpage
\section{Introduction}

\begin{quotation}
{\footnotesize \it
`This track, as you perceive, was made by a rider who was going from the direction of the school.' 

`Or towards it?' 

`No, no, my dear Watson. The more deeply sunk impression is, of
course, the hind wheel, upon which the weight rests. You perceive
several places where it has passed across and obliterated the more
shallow mark of the front one. It was undoubtedly heading away from
the school.' } \small

\vspace{.1in} \noindent From {\em The Adventure of the Priory School}, a 
Sherlock Holmes story by A. Conan Doyle.

\end{quotation}

 \noindent
The  radial projection
 of a Brownian motion started at the origin and run for 
 unit time in $d$ dimensions 
defines a random occupation measure on the sphere $S^{d-1}$.
 Can we determine the
endpoint of the Brownian path from this projected occupation measure?
The problem of recovering data given a projection of the data is a
common theme both inside and outside of probability theory.
The title of this paper is adapted from a handout distributed
 by Peter Doyle, where 
the geometric problem of recovering from bicycle tracks the exit
direction of the cyclist was posed.  

An interesting feature of the present reconstruction problem is that
the answer in low dimensions is different from the answer in dimensions
$d \geq 3$.  This would not be too surprising, except that the
behavior in the one-dimensional case involves a conditioning
identity which does not seem inherently one-dimensional.  This identity
concerns the
conditional distribution of the endpoint given the occupation measure.
One of the aims of this paper is to understand why this identity
breaks down in higher dimensions,
and what version of this identity
might hold even when the occupation measure determines the endpoint and
indeed determines the entire unprojected path.  In high dimensions,
recovery of the endpoint (and entire path), while intuitively plausible,
is somewhat tricky because, as described in \cite[page 275]{IM74}, 
the particle ``comes in spinning''.  In particular, the
range of the projected path is a.s. a dense subset of the sphere
(see remark at the end of this introduction).
Thus some quantitative criterion on accumulation of measure is required
even to recover the set of occupied points on the sphere from the occupation 
measure.

Throughout the paper $d$ is a positive integer, and
$S^{d-1} \subseteq \Rd$ is the unit sphere.  
We often omit $d$ in the notation for various spaces and 
mappings whose definition depends on $d$. 
Let $\sph : \Rd \to S^{d-1}$ be the spherical
projection $\sph (x) = x / |x|$ for $x \ne 0$, with some arbitrary conventional
value for $\sph(0)$.
Let $(W_t, t \ge 0)$ denote a standard Brownian motion in $\Rd$ with $W_0 = 0$,
which we take to be defined on some underlying probability space 
$(\Omega , \F , \P)$.
For $t \ge 0$ let $\Theta_t:= \sph(W_t)$, and let
$\Theta:= (\Theta_t, 0 < t \le 1 )$.
Let $\mu_\Theta$ denote the random occupation measure of $\Theta$ on
$S^{d-1}$, that is
\eq
\lb{defmu}
\mu_\Theta (B):= \int_0^1 \one_{\Theta_t \in B} \, dt
\en
for Borel subsets $B$ of $S^{d-1}$. We may regard $\mu_\Theta$ as a
random variable defined on $(\Omega , \F , \P)$ with values in the space 
$(\omxxx , \F_2)$ of Borel probability measure on
$S^{d-1}$ endowed with the $\sigma$-field generated by the measures
of Borel sets.  

The questions considered in this paper arose from the following identity:
for each Borel subset $B$ of $S^{d-1}$, we have
\eq
\lb{cform}
P( \Theta_1 \in B  \,|\,\mu_\Theta (B) ) = \mu_\Theta (B).
\en
If $d =1$ then $S^{0} = \{-1,1\}$, and
$\mu_\Theta (\{1\})$ and $\mu_\Theta (\{-1\})$ are the times spent positive
and negative respectively by a one-dimensional Brownian motion up to
time 1. As observed in Pitman-Yor \cite{py92}, formula \re{cform} in this case
can be read from \Lev's description \cite{lev39} of the joint law of the arcsine distributed 
variable $\mu_\Theta (\{1\})$ and the Bernoulli$(1/2)$
distributed indicator $\one_{W_1 >0}$.
The truth of \re{cform} in higher dimensions is not so easily checked,
due to the lack of explicit formulae for the distribution of
$\mu_\Theta (B)$ even for the simplest subsets $B$ of $S^{d-1}$;
see for instance \cite{bd88}.
However, \re{cform} can be deduced from the scaling
property of Brownian motion, which implies that the process
$(\Theta_t, t \ge 0)$ is {\em $0$-self-similar}, meaning 
there is the equality in distribution
$$(\Theta_t, t \ge 0 ) \ed (\Theta_{ct}, t \ge 0 )$$
for all $c >0$.
According to an identity of Pitman and Yor \cite{py95rdd},
recalled as Proposition \ref{pr:PY} in Section \ref{ids}, 
the identityi \re{cform} holds for an arbitrary jointly measurable 
$0$-self-similar process $(\Theta_t, t \ge 0)$ with values in an 
abstract measurable space, for any measurable subset $B$ of that space.

Formula \re{cform} led us to the following question, which we discuss
further in Section \ref{quest}:
\begin{qst}{q1}
For which processes $\Theta:= (\Theta_t, 0 \le t \le 1)$, does the
identity
\eq
\lb{cform1}
P( \Theta_1 \in B  \,|\,  \mu_\Theta) = \mu_\Theta (B)
\en
hold for all measurable subsets $B$ of the range space of $\Theta$? 
\end{qst}
To clarify the difference between \re{cform} and \re{cform1},
$P( \Theta_1 \in B  \,|\,  \mu_\Theta(B))$ on
the LHS of \re{cform1} is a conditional probability given 
the $\sigma$-field generated by the real random variable $\mu_\Theta(B)$,
whereas $P( \Theta_1 \in B  \giv  \mu_\Theta)$ on
the LHS of \re{cform1} is a conditional probability given the 
$\sigma$-field generated by the random measure $\mu_\Theta$, that is
by all the random variables $\mu_\Theta(C)$ as $C$ ranges over measurable subsets of the range space of $\Theta$.
For a general process $\Theta$, formula \re{cform1} implies \re{cform},
but not conversely.

Now let $\Theta$ be the spherical projection of Brownian motion.
If $d =1$ then the $\sigma$-field generated by $\mu_\Theta$
is identical to that generated by either $\mu_\Theta (\{1\})$ or by
$\mu_\Theta (\{-1\}) = 1 - \mu_\Theta (\{-1\})$. So \re{cform1} is
a consequence of \re{cform} if $d=1$. 
But \re{cform1} fails for $d \ge 2$. We show this for $d=2$
in Section \ref{quads} by some explicit estimates involving the occupation
times of quadrants. For $d \ge 3$ formula \re{cform1} fails even more
dramatically.  In Section \ref{recov} we show that if $d \ge 3$ then
$\Theta_1$ is a.s. equal to a measurable function of $\mu_\Theta$.
Less formally, we say that $\Theta_1$ {\em can be recovered from $\mu_\Theta$}.
This brings us to the question of what features 
of the path of the original Brownian motion $W:= (W_t, 0 \le t \le 1)$ can be
recovered from $\mu_\Theta$.
If $d =3$ it is well known that
the path $W$ has self-intersections almost surely,
so one can define a measure-preserving map $T$ on the Brownian path space that 
reverses the direction of an appropriately selected closed loop in the path.
Regarding $\mu_\Theta = \mu_\Theta(W)$ as a function on path space, we then
have $\mu _\Theta (W) = \mu_\Theta  (T W)$, hence
$\P (W \in A \| \mu_\Theta) = \P (W \in T^{-1} A \| \mu_\Theta)$,
from which it follows that $W$ itself cannot be recovered from $\mu_\Theta$.
However, for $d =3$ it is possible to recover 
from $\mu_\Theta$ both the random set
$$
\range (W) := \{ W_t : 0 \leq t \leq 1 \}
$$
and the final value $W_1$. 
We regard $\range (W)$ as a random variable with values in the space
$\omxxxx$ of closed bounded subsets of $\Rd$, 
equipped with the Borel $\sigma$-field for the Hausdorff metric 
$$d(S,T) := \max \{ \sup_{x \in S} d(x,T) \, , \, \sup_{y \in T} d(y,S) \}$$
where $d(x,S) := \inf_{y \in S} |x-y|$. 
We now make a formal statement of the recovery result:

\begin{thm} {th:recover}  Fix the dimension $d \geq 3$. 
There exist measurable functions $$\psi : \omxxx \to \omxxxx \; \mbox{ \rm and } \:
\varphi : \omxxxx \to \Rd ,$$
 such that there are almost sure
equalities
\begin{eqnarray} \label{eq:psi}
\psi (\mu_\Theta ) & = & \range (W) \\
\mbox{and } \hspace{.6in} \varphi \circ \psi (\mu_\Theta ) & = & 
   W_1. \label{eq:phi}
\end{eqnarray}
\end{thm}

For $d \ge 4$ it is a routine consequence
of the almost sure parameterizability of a Brownian path by its quadratic
variation that $W$ can be recovered from $\range (W)$.
So Theorem \ref{th:recover} implies that $W$ can be recovered from
$\mu_\Theta$ in dimensions $d \ge 4$.
The only part of the proof of Theorem \ref{th:recover} 
which involves probabilistic estimates for Brownian motion is Lemma~\ref{lem:rho}, the remainder being mostly 
point-set topology.  We remark that the topological arguments below
also show that in dimension $d=2$, the range of $W$ and 
the endpoint $W_1$ can be recovered 
from the occupation measure of the planar path $W:= (W_t, 0 \le t \le 1)$

As usual, the hardest (and most interesting) dimension is two.  We
conjecture that the high-dimensional behavior does not extend
to two dimensions, that is,

\begin{cnj} {th:no recover}
When $d = 2$, there is no map $\psi : \omxxx \to \omxxxx$ such that almost
surely 
$$\psi ( \mu_\Theta )  = \range (W) .$$
\end{cnj}

When $d=2$ it can be deduced from work of Bass and Khoshnevisan \cite[Theorem 2.9]{bassk92} that $\mu_\Theta$ almost surely has a continuous density, call it
the {\em angular local time process}. The problem of describing the conditional
law of $W$ given $\mu_\Theta$ for $d=2$ is then analogous to the problem
studied by Warren and Yor \cite{wy98}, who give an account of the 
randomness left in a one-dimensional Brownian motion after conditioning on its occupation measure up to a suitable random time.
Aldous~\cite{al97} and Knight~\cite{kn98} treat related
questions involving the distribution of Brownian motion conditioned on
 its local time process.
However, as far as we know there is no Ray-Knight type description 
available for the angular local time process, and this makes it
difficult to settle the conjecture.

\noindent{\bf Remark.} Let $\Theta_t:=W_t/|W_t|$ be the radial
projection of Brownian motion in $\R^d$. It is a classical fact
that for any $\epsilon>0$, the initial path segment 
$\{\Theta_t \, : \, 0<t<\epsilon\}$ is dense in the unit sphere $S^{d-1}$.
Since this fact motivates much of our work, we include an
elementary explanation for it, which is valid in greater generality.
It suffices to show that for any open set $U$ on the sphere
and any $\epsilon>0$,
the probability of the event $E(U,\epsilon)$ that
$\{\Theta_t \, : \, 0<t<\epsilon\}$ intersects $U$, equals one.
By compactness, some finite number $N_U$ of rotated copies of $U$
cover the sphere, so by rotation invariance of Brownian motion,
$\P[E(U,\epsilon)]\ge N_U^{-1}$. Therefore
$$\P\Big[\bigcap_{n=1}^\infty E(U,1/n)\Big]\ge N_U^{-1} \,,
$$
whence by the Blumenthal zero-one law, this probability must be 1.

\section{Identities for scalar self-similar processes}
\label{ids}
\setcounter{equation}{0}

Recall that a real or vector-valued process $(X_t,t \ge 0)$ is
called {\em $\beta$-self-similar} for a $\beta \in \reals$
if for every $c > 0$
\eq
\lb{selfsimdef}
( X_{ct}, t \ge 0 ) \ed ( c^{\beta} X_{t}, t \ge 0 )
\en
Such processes were studied by Lamperti \cite{lam62ss,lam72zw},
who called them {\em semi-stable}.
See \cite{taq86} for a survey of the literature of these processes.
The conditioning formula \re{cform} for any $0$-self-similar process
$(\Theta_t, t \ge 0)$ is an immediate consequence of the following 
identity. 
To see the direct implication,
take $X(t , \omega)$ to equal $\one_{(0,\infty)} (\omega (t))$.

\begin{prp} {pr:PY} {\em (Pitman and Yor \cite{py95rdd})}
Let $( X_t , t \geq 0 )$ be stochastic process with
$X : \R^+ \times \Omega \to \R$ jointly measurable.  
Let $\Xb_t:= t^{-1} \int_0^t X_s \, ds$ and suppose that
\begin{eqnarray}
(X_t , \Xb_t) & \diseq & (X_1 , \Xb_1) \label{eq:self similar} \\
\mbox{and   } \;\;\; \E |X_1| & < & \infty \label{eq:integrable} .
\end{eqnarray}
Then for every $t > 0$,
\begin{equation} \label{eq:ss conclusion}
\E (X_t \| \Xb_t) = \Xb_t .
\end{equation}
\end{prp}

\noindent{\sc Proof:}  We simplify slightly the proof
in \cite{py95rdd}.  Due to~(\ref{eq:self similar}) it suffices to
prove~(\ref{eq:ss conclusion}) for $t = 1$.  It also suffices to
prove this on the event $\{ \Xb_1 \neq 0 \}$, since this implies 
$\E X_1 \one_{\Xb_1 \neq 0} = \E \Xb_1 \one_{\Xb_1 \neq 0}$ and
subtracting the relation $\E X_1 = \E \Xb_1$ (a consequence
of~(\ref{eq:self similar})) yields $\E X_1 \one_{\Xb_1 = 0} = 0$.  
This is equivalent to proving
\begin{equation} \label{eq:f version}
\E [f(\Xb_1) ; \Xb_1 \neq 0 ] = \E \left [ f(\Xb_1) {X_1 \over \Xb_1} ; 
   \Xb_1 \neq 0 \right ]
\end{equation}
for a suitably large class of functions $f$.  Let $\nu$ be the law of
$\Xb_1$.  Since $f(x) \one_{x \neq 0}$ for bounded measurable $f$ 
may be approximated in $L^2 (\nu)$
by bounded functions vanishing in a neighborhood of zero and having bounded 
continuous derivative, this class suffices.  Fix such a function $f$
and apply the chain rule for Lebesgue integrals (see, e.g., \cite{RY94}, 
Chapter 0, Prop.~(4.6)), treating $\omega$ as fixed, to obtain
$$f \left(\int_0^1 X_t \, dt\right) = \int_0^1 f' \left(\int_0^t X_s \, ds\right) X_t \, dt .$$
Boundedness of $f'$ allows the interchange of expectation with integration,
so using~(\ref{eq:self similar}) we get
(\ref{eq:f version}) from the following computation:
\begin{eqnarray*}
\E [f(\Xb_1) ; \Xb_1 \neq 0 ] = 
\E f(\Xb_1) & = & \int_0^1 \E \left [ f' \left(\int_0^t X_s \, ds\right) X_t 
   \right] \, dt \\
& = & \int_0^1 \E \left[ f' (t \Xb_1)  X_1 \right] \, dt \\
& = & \E \left[ \int_0^1 f' (t \Xb_1)  X_1 \, dt \right] \\
&=&\E \left[ f(\Xb_1) {X_1 \over \Xb_1} \right]\, . 
\end{eqnarray*}
\endpf

For a different proof and variations of the identity see \cite{py94mn60}.
We see immediately that~(\ref{eq:self similar}) holds for any 
0-self-similar process $X$.  We observe also:

\begin{crl} {cor:cones}
Let $( Y_t )$ be any $\beta$-self-similar vector-valued process.
Let $X_t := \one_{\{Y_t \in \C\}}$ where $\C$ is any Borel set which is a
cone, i.e., for $\lambda > 0$, $x \in \C \Leftrightarrow \lambda x \in \C$.  
Then 
$( X_t )$  satisfies~{\em (\ref{eq:self similar})}, and hence
\begin{equation} \label{eq:cone}
\P (Y_1 \in \C \| \Xb_1) = \Xb_1 .
\end{equation}
\end{crl}


Applying Bayes' rule to~(\ref{eq:cone}) yields the following corollary.

\begin{crl} {cor:Bayes}
Let $\{ Y_t \}$  be any $\beta$-self-similar vector-valued process,
and let $V_t = \int_0^t X_s \, ds$ with $X_t := \one_{\{Y_t \in \C\}}$
for a fixed positive cone $\C$.  Then 
$$\P (V_t \in dv \| Y_t \in \C) = {v \P (V_t \in dv) \over 
   t \P (Y_t \in \C)} \, .$$
\end{crl}

\begin{crl} {cor:beta}
Under the hypotheses of the Corollary \ref{cor:Bayes}, suppose $\Xb_t$ has
a ${\rm beta}(a,b)$ distribution.  Then the conditional distribution of
$\Xb_t$ given $Y_t \in C$ is ${\rm beta}(a+1 , b)$ and the conditional
distribution of $\Xb_t$ given $Y_t \notin \C$
 is ${\rm beta}(a , b+1)$.  
\end{crl}

\begin{exm}{Levy} Stable L\'evy Processes.
{\em
Let $\{ Y_t \}$ be a stable L\'evy process
that satisfies $\, \P (Y_t > 0) = p$ for all $t$.  It is well known \cite{pr69,gs94}
that the distribution of the total duration $V_1$ 
that $\{ Y_t \}$ is positive up to time 1, is ${\rm beta}(p , 1-p)$.
It follows that the conditional distributions of $V_1$ given
the sign of $Y_1$ are respectively
\begin{eqnarray}
(V_1 \| Y_1 > 0) & \diseq & {\rm beta}(1 + p , 1 - p) \label{eq:beta1} \\
(V_1 \| Y_1 < 0) & \diseq & {\rm beta}(p , 2 - p)  \label{eq:beta2} .
\end{eqnarray}
}
\end{exm}

\begin{exm}{perturbed} Perturbed Brownian Motions.
{\em
Let $Y_t:= |B_t | - \mu \ell_t, t \ge 0 $, where $B$ is a standard one-dimensional
Brownian motion started at $0$, $\mu > 0 $ and 
$(\ell_t, t \ge 0 )$ is the local time process of $B$ at zero.
F. Petit \cite{pet92a} showed that $V_1^- := \int_0^1 ds 1 _{ (Y_s < 0)}$
has beta$\left( {1 \over 2 } ,{ 1 \over 2 \mu } \right)$ distribution.
Corollary \ref{cor:beta} implies that 
the conditional distribution of $V_1^-$ given $Y_1 < 0$ is
beta$\left( {3 \over 2 } ,{ 1 \over 2 \mu } \right)$ 
and that
the conditional distribution of $V_1^-$ given $Y_1 > 0$ is
beta$\left( {1 \over 2 } , 1 +{ 1 \over 2 \mu } \right)$. 
These results have been stated and proved in
\cite[Th. 8.3]{yor92z} and in \cite{cpy93}.
A more general class of beta laws has been obtained for the times spent
in $\reals_\pm$ by doubly perturbed Brownian motions, that is to say solutions
of the stochastic equation
$$
Y_t = B_t +  \alpha \sup_{0 \le s \le t}Y_s +  \beta \inf_{0 \le s \le t}Y_s  .
$$
See, {e.g.}, Carmona-Petit-Yor \cite{cpy96},
Perman-Werner \cite{perwer96} and
Chaumont-Doney \cite{chdon96}.
}
\end{exm}

\begin{exm}{levmore} More about the Brownian case.
{\em
Formula \re{cform}
has some surprising consequences even in the simplest case when
$d=1$.
Consider the function
\eq
\lb{fta}
f(t,a) := P(B_t >0 | V_1 = a)
\en
for $0 < t \le 1 $ and $0 \le a \le 1$, where $B$ is a one-dimensional
Brownian motion and
$V_1 = \int_0^1 1(B_t > 0 ) dt$.
Without attempting to compute $f(t,a)$ explicitly, 
which appears to be quite difficult, let us
presume that $f$ can be chosen to be continuous in $(t,a)$.
Then
\eq
\int_0^1 f(t,a) dt = a = f(1,a )  ~~~~~~~~~~~(0 \le a \le 1 )
\en
where the first equality follows from \re{fta} and the second equality
is read from \re{cform}.
On the other hand, it is easily seen that
\eq
f(0+,a ) = \hf               ~~~~~~~~~~~~~~~~(0 < a < 1 )
\en
which implies that 
\eq
\mbox{ \em for each $a > \hf $ there exists $t \in (0,1)$ such that  
$f(t,a) > a$ }
\en
That is to say, given $V_1 = a > \hf$, there is some time $t < 1$ such that
the BM is more likely to be positive at time $t$ than it is at time $1$.
}
\end{exm}

\section{Identities for self-similar processes in dimension $d \ge 2$}
\setcounter{equation}{0}
\label{quest}

Say that a jointly measurable process $\Theta:= (\Theta_t, 0 < t \le 1 )$ has the
{\em sampling property} if
\eq
\lb{cform11}
\P( \Theta_1 \in B  \,|\,  \mu_\Theta) = \mu_\Theta (B)
\en
for all measurable subsets $B$ of the range space of $\Theta$.
The results of this section consist
of two examples where the sampling property does hold, and a
characterization of the sampling property in terms of exchangeability.

\begin{prp} {pr:exchangeable}
Suppose that $\Theta$ takes values in a Borel space.
Let $U_1, U_2, \cdots $ be a sequence
of i.i.d. random variables with uniform distribution on $(0,1)$,
independent of $\Theta$.
Then the following are equivalent:

(i) $(\Theta_t)$ has the sampling property;

(ii) for each $n = 1, 2, 3, \cdots $, 
\begin{equation} \label{eq:exchangeable}
(\Theta_1, \Theta_{U_2}, \Theta_{U_3}, \cdots, \Theta_{U_n} )
\ed (\Theta_{U_1}, \Theta_{U_2}, \Theta_{U_3}, \cdots, \Theta_{U_n} ) .
\end{equation} 
\end{prp}

\noindent{\sc Proof:} 
Clearly (ii) is equivalent to
\eq
\lb{equiv}
\P( \Theta_1 \in B  \,|\, \{\Theta_{U_j}\}_{j=2}^\infty ) =
\P( \Theta_{U_1} \in B  \,|\, \{\Theta_{U_j}\}_{j=2}^\infty )
\en
for all measurable subsets $B$ of the range space of $\Theta$.
To connect this to (i), observe that 
$\{\Theta_{U_j}\}_{j=2}^\infty$ is a sequence 
of i.i.d.\ picks from $\mu_\Theta$. Hence this sequence
is conditionally independent of $\Theta$
given $\mu_\Theta$.
Therefore, (\ref{equiv})
can be rewritten as 
\eq
\lb{equiv2}
\P( \Theta_1 \in B  \,|\,\mu_\Theta ) =
\P( \Theta_{U_1} \in B  \,|\, \mu_\Theta )
\en
for all measurable $B$, which is equivalent to (i).
\endpf

The conditions \re{eq:exchangeable} increase in strength as $n$ increases.  
For $n=2$,~(\ref{eq:exchangeable}) is just 
\begin{equation} \label{eq:n=2}
(\Theta_1 , \Theta_{U_2}) \diseq (\Theta_{U_1} , \Theta_{U_2}) .
\end{equation}
which immediately implies
\begin{equation} \label{eq:two-exch}
(\Theta_1 , \Theta_{U_2}) \diseq (\Theta_{U_2} , \Theta_1) .
\end{equation}

\begin{prp} {pr:two-exch}
If the distribution of $(\Theta_s , \Theta_t)$ depends only on $t/s$ then 
the conditions~{\em (\ref{eq:n=2})} and~{\em(\ref{eq:two-exch})} are equivalent.
\end{prp}

\noindent{\sc Proof:} Construct $U_1$ and $U_2$ as follows.
Let $Y$ and $Z$ be independent with $Y$ uniform on $[0,1]$
and $Z$ having density $2x$ on $[0,1]$.  Let $X$ be an independent
$\pm 1$ fair coin-flip and set $(U_1 , U_2)$ equal to
$(Z , YZ)$ if $X = 1$ and $(YZ , Z)$ if $X = -1$.  By construction,
the law of $(\Theta_{U_1} , \Theta_{U_2})$ is one half the law of 
$(\Theta_Z , \Theta_{YZ})$ plus one half the law of 
$(\Theta_{YZ} , \Theta_Z)$.  By the assumption on $\Theta$
this is one half the law of $(\Theta_1 , \Theta_{U_2})$ plus one half
the law of $(\Theta_{U_2} , \Theta_1)$.  This and~(\ref{eq:two-exch})
imply~(\ref{eq:n=2}).  
\endpf

We note that the spherical projection of Brownian motion in $\Rd$ 
satisfies (\ref{eq:two-exch}) for all $d$.
So this condition is not enough to imply the sampling property for
a $0$-self-similar process $\Theta$.
%
%
When $\Theta$ is not 0-self-similar it is easy to find
cases where~(\ref{eq:two-exch}) holds but not~(\ref{eq:n=2}).  
\begin{exm} {converse}
{\em
Let $(X,Y)$ have a symmetric distribution and let
$\Theta_t = X \one_{t < a} + Y \one_{t \geq a}$ for a fixed $a \in (0,1)$.
It is easy to see that~(\ref{eq:two-exch}) holds.  On the other hand,
if $\P (X = Y) = 0$, then $\P (\Theta_1 = \Theta_{U_2}) = 1-a$ while
$\P(\Theta_{U_1} = \Theta_{U_2}) = a^2 + (1-a)^2$.  Unless $a = \hf$,
these two probabilities are not equal.  
}
\end{exm}

We now mention some interesting examples of $0$-self-similar processes 
which do have the sampling property.

\begin{exm}{walsh}
{Walsh's Brownian motions.}
{\em
Let $B$ be a one-dimensional BM started at 0. Suppose that each excursion
of $B$ away from $0$ is assigned a random angle in $[0,2 \pi)$ according to 
some arbitrary distribution, independently of all other excursions.
Let $\Theta_t$ be the angle assigned to the excursion in progress at time
$t$, with the convention $\Theta_t = 0$ if $B_t = 0$.
So $( |B_t|, \Theta_t)$ is Walsh's singular Brownian motion in the plane 
\cite{wal78,bpy89a}.
As shown in \cite[Section 4]{py92},
the process $(\Theta_t)$ is a $0$-self-similar process with the 
sampling property, and the same is true
of $(\Theta_t)$ defined similarly for a $\delta$-dimensional Bessel process
instead of $|B|$ for arbitrary $0 < \delta < 2$.
}
\end{exm}
The proof of the sampling property of the angular part $(\Theta_t)$ of 
Walsh's Brownian motion is based on the following lemma, which is implicit
in arguments of \cite[Section 4]{py92} and \cite[formula (24)]{py95rdd}.

\begin{lmm} {lmmclosed}
Let $Z$ be a random closed subset of $[0,1]$ with Lebesgue measure zero. 
For $0 \le t \le 1$
let $N_t - 1$ be the number of component intervals of the set
$[0,t] \backslash Z$ whose length exceeds $t - G_t$, where $G_t = \sup \{s : s < t, s \in Z \}.$
So $N_t$ has values in $\{1,2, \cdots , \infty \}$.
Given $Z$,
let $(\Theta_t)$ be a process constructed by assigning each complementary
interval of $Z$ an independent angle
according to some arbitrary distribution on $[0, 2\pi)$, and letting
$\Theta_t = 0$ if $t \in Z$.
If $(N_t)$ has the sampling property, then so does $(\Theta_t)$.
\end{lmm}

According to \cite[Theorem 1.2]{py92} and
\cite[formula (24)]{py95rdd},
for $Z$ the zero set of a Brownian motion,
or more generally the range of a stable$(\alpha)$ subordinator for $0 < \alpha <1$, the process $(N_t)$ has the sampling property, hence so does
the angular part $(\Theta_t)$ of Walsh's Brownian motion whose radial
part is a Bessel process of dimension $\delta$ for arbitrary $0 < \delta < 2$.

\begin{exm}{dirich}
{A Dirichlet Distribution.}
{\em
Let $Z$ be the set of points of a Poisson random measure
on $(0,\infty)$ with intensity measure $\theta x^{-1} dx, x > 0$.
Construct $(\Theta_t)$ from $Z$ as in Lemma \ref{lmmclosed}.
So between each pair of points of the Poisson process, an independent
angle is assigned, with some common distribution $H$ of angles on $[0, 2 \pi)$.
It was shown in \cite{py95rdd} that $(N_t)$ derived from this $Z$
has the sampling property,
hence so does $(\Theta_t)$ derived from this $Z$.
In this example $\mu_\Theta$ 
is a Dirichlet random measure governed by $\theta H$
as studied in \cite{fe74,ki75,ig82,seth94}.
}
\end{exm}

We close this section by rewriting Proposition~\ref{pr:PY} 
as a statement concerning stationary processes.
Let $( X_t )$ be a
jointly measurable process and $Y_t = X_{e^t}$.  The process
$( X_t )$ being 0-self-similar is equivalent to the process
$( Y_t )$ being stationary, so a change of variables
turns Proposition~\ref{pr:PY} into:
\begin{crl} {cor:stationary}(Pitman-Yor \cite{py94mn60})
Fix $\lambda > 0$ and define $\Yb_\lambda := \int_0^\infty \lambda 
e^{-\lambda t} Y_t \, dt$, where $\{ Y_t : t \in \R \}$ is a 
stationary process and $\E |Y_0| < \infty$.  Then
$$\E (Y_0 \| \Yb_\lambda) = \Yb_\lambda .$$
\end{crl}

The following proposition provides a partial converse:

\begin{prp} {pr:characterize} 
Let $F$ be a distribution on $[0 , \infty)$ and for a stationary
process $( Y_t )$ let $\Yb_F$ denote $\int_0^\infty Y_t \, dF$.
Assuming either $F$ has a density or $F$ is a lattice distribution,
the identity $\E (Y_0 \| \Yb_F) = \Yb_F$ holds for every such process
$\{ Y_t \}$ if and only if $F$ has density $\lambda e^{-\lambda t}$
for some $\lambda \in (0 , \infty )$ or $F = \delta_0$.
\end{prp}

\noindent{\sc Proof:} Fix $F$ and suppose that $\E (Y_0 \| \Yb_F) = \Yb_F$
holds for all stationary $\{ Y_t \}$ with $\E |Y_0| < \infty$.
When also $\E |Y_0|^2 < \infty$, this implies $\E Y_0 \Yb_F = \E (\Yb_F)^2$.
Let $r(t) = \E Y_0 Y_t$ and let $\xi_1 , \xi_2$ be i.i.d. according to $F$.  
Comparing 
$$\E Y_0 \Yb_F = \E r(|\xi_1|)$$
with
$$\E (\Yb_F)^2 = \E r(|\xi_1 - \xi_2|) , $$
we find that 
$$\E r(|\xi_1|) = \E r(|\xi_1 - \xi_2|) .$$
Taking $Y$ to be an Ornstein-Uhlenbeck process shows that this holds
for $r(t) = e^{-\alpha t}$, so that $|\xi_1 - \xi_2|$ has the
same Laplace transform, hence the same distribution, as $|\xi_1|$.
Assuming that $F$ is concentrated on $[0,\infty)$ and has a density,
Puri and Rubin \cite{purirubin70} showed that this condition implies
$F$ is an exponential.  If $F$ is a lattice distribution, they
showed it must be $\delta_0$ or $\hf \delta_0 + \hf \delta_a$ or
$a$ times a geometric for some $a > 0$.  It is easy to construct
examples ruling out the nondegenerate discrete cases.  
\endpf

Changing back to $X_t := Y_{\log t}$, Proposition \ref{pr:characterize} yields:

\begin{crl} {cor:ss}
Suppose $F$ has a density
$f$ on $(0,1)$. 
The identity $$\E (X_1 \| \int_0^1 X_s \, dF) = \int_0^1 X_s \, dF$$
holds for all 0-self-similar processes $( X_t )$ with $\E |X_1| < \infty$
if and only if $f(x) = \lambda x^{\lambda - 1}$ for some $\lambda > 0$.
\end{crl}

\section{Quadrants and the two-dimensional case}
\label{quads}
\setcounter{equation}{0}

In this section we establish the following
Proposition.

\begin{prp} {pr:2 dim}
Let $d = 2$ and let $Q_1 , Q_2 , Q_3 , Q_4$ be the four quadrants
in the plane, in clockwise order.  
Let 
$$
\mu(Q_i):= \int_0^1 \one_{W_t \in Q_i} \, dt
$$
denote the time spent in $Q_i$ up to time 1 by a planar Brownian motion
$W$ started at the origin.
Then for each $k \leq 4$, the random variable
$$\P (W_1 \in Q_k \| \mu (Q_i) : 1 \leq i \leq 4) $$
is not almost surely equal to $\mu (Q_k)$.  
\end{prp}

Fix the
dimension $d=2$ throughout, and denote by $A_{\ee}$ the event that
$\mu(Q_2) \in [\ee, 2\ee]$,
$\mu(Q_3) \in [\ee, 2\ee]$, and
$\mu(Q_4) \in [\ee, 2\ee]$. 
Thus, if
$A_{\ee}$ occurs, then the Brownian motion $W$ spends only a small
amount of time in $Q_2, Q_3$, and $Q_4$.  The idea behind the
proof is that if Brownian motion spends most of its time in $Q_1$, then it
is very unlikely to be in $Q_3$ at time $1$, since $Q_1$ and $Q_3$
do not share a common boundary.  More precisely, we will show that
there is a constant $C$ for which 
\begin{equation} \label{eq:quadrants}
\P (W_1 \in Q_3 | A_{\ee}) \leq C \ee^2 [\log(1/\ee)]^3
\end{equation}
for sufficiently small $\ee > 0$, which clearly 
implies  Proposition~\ref{pr:2 dim}. 
The estimate \re{eq:quadrants} follows immediately from
the lower bound for $\P (A_{\ee})$ and the upper bound for
$\P (\{ W_1 \in Q_3 \} \cap A_{\ee})$
given in Lemmas \ref{thirdlemma} and \ref{fourthlemma} below.

\begin{lmm} {firstlemma}
Let $( B_t )$ be one-dimensional Brownian motion started from the
origin.  Then as $\delta \te 0$
\eq
\lb{con1}
\delta^{-1} \P ( \min_{t \in [0,1]} B_t \geq - \delta)
   \rightarrow \sqrt{2 \over \pi}
\en
and
\eq
\lb{con2}
\delta^{-3} \P (\min_{t \in [0,1]} B_t \geq - \delta \mbox{ and } B_1 < 0) 
\rightarrow {1 \over \sqrt{2 \pi}}
\en
\end{lmm}

\noindent{\sc Proof}:  The first limit results
from the fact that $\min_{t \in [0,1]} B_t$ has density $2 \phi(x)$ on
$(-\infty,0]$ where $\phi$ is the standard normal density of $B_1$.
The second follows easily from the reflection principle, which shows that the 
probability involved equals
$$\int_{-\delta}^0 (\phi (x) - \phi (x - \delta) )\, dx$$.

\endpf

\begin{lmm} {thirdlemma} 
There exists a constant $C_1 > 0$ such that
$$\P (A_{\ee}) \geq C_1 \ee$$ 
for sufficiently small $\ee > 0$.
\end{lmm}

\noindent{\sc Proof}: Let $D_\ee$ be the set $(\sqrt{\ee} , \sqrt{\ee}) + Q_1$
and let $S_{\ee}$ be the event that
$$\{ \left | \{ t \in [0, 4\ee]: W_t \in Q_i \} \right | \geq \ee
   \mbox{ for } i = 2,3,4, \mbox{ and } W_{4 \ee} \in D_\ee \} .$$
Let $C_2 = \P (S_{\ee}) > 0$.  From the scaling properties of
Brownian motion, we see that $C_2$ does not depend on $\ee$.
Let $p_\ee$ be the probability that $\min_{t \in [0 , 1 - 4\ee]} B_t 
> - \sqrt{\ee}$.  By the Markov property and independence of the
coordinates of $W$, $\P (A_\ee) \geq C_2 p_\ee^2$.  
Lemma~\ref{firstlemma} tells us that $p_\ee \geq \sqrt{ (2 - \beta) \ee / \pi}$
for any $\beta > 0$ and sufficiently small $\ee (\beta)$.  This  proves 
the lemma with any $C_1 < 2 C_2 / \pi$.   
\endpf

\begin{lemma}
\label{fourthlemma} There exists a constant $C_3 < \infty$ such that
$$\P (\{ W_1 \in Q_3 \} \cap A_{\ee}) \leq C_3 \ee^3 [\log(1/\ee)]^3$$ 
for sufficiently small $\ee > 0$.
\end{lemma}

\noindent{\sc Proof}:  Choose $C_4 > 12$ and let
$\delta = C_4 \sqrt{\ee \log(1/\ee)}$.  Let $Q_1^{\delta} =
\{ (x,y): x > -\delta, y > -\delta \}$.  Also define $T_{\delta} =
\min \{t: W_t \notin Q_1^{\delta} \}$.  Let $R_1 = 
A_{\ee} \cap \{ T_{\delta} \leq 1 - 6\ee \}$, $R_2 =
A_{\ee} \cap \{1 - 6\ee < T_{\delta} \leq 1 \}$, and
$R_3 = \{ W_1 \in Q_3 \} \cap \{ T_{\delta} > 1 \}$.  By splitting
up the event $\{ W_1 \in Q_3 \} \cap A_{\ee}$ according to
the value of $T_{\delta}$, we see that
if $\{ W_1 \in Q_3 \} \cap A_{\ee}$ occurs, then either
$R_1$, $R_2$, or $R_3$ must occur.  We will prove the lemma by
establishing upper bounds on $\P (R_1)$, $\P (R_2)$, and
$\P (R_3)$.

To bound $\P (R_3)$, apply \re{con2} to the two
independent coordinate processes, yielding for sufficiently
small $\ee$
$$\P (R_3) \leq \delta^6 = {C_4^6 } \ee^3 \log (1/\ee)^3. $$

A bound for $\P (R_2)$ follows from the observation that on $A_\ee$,
there must be some $t \in [1 - 6 \ee , 1]$ for which $W_t \in 
Q_1$.  Thus on $R_2$, one of the two coordinate processes has an 
oscillation of at least $\delta$ on the time interval $[1 - 6 \ee , 1]$.
This implies that one of the coordinate processes strays by at least
$\delta / 2$ from its starting value in the interval $[1 - 6 \ee , 1]$,
hence by the Markov property,
\begin{eqnarray*}
\P (R_2) & \leq & 2 \P (\max_{0 \leq t \leq 6 \ee} |B_t| \geq {\delta / 2})
   \\[2ex]
& \leq & 8 \P ( B_{6 \ee} \geq {\delta / 2}) , \mbox{ by the reflection
principle,} \\[2ex]
& \leq & 8  \exp (- \delta^2 / 48 \ee) = 8 \ee^{C_4^2/48} .
\end{eqnarray*}
By choice of $C_4 >12$, this is $o(\ee^3 )$.

A bound on $\P (R_1)$ may be obtained in a similar way.  Observe that
on $A_\ee$, there must be some $t \in [T_\delta , T_\delta + 6 \ee]$
for which $W_t \in Q_1$.  Thus one of the coordinates increases
by at least $\delta$ from its starting value on the time interval 
$[T_\delta , T_\delta + 6 \ee]$.  The strong Markov property yields
$$\P (R_1) \leq 2 \P (\max_{0 \leq t \leq 6 \ee} B_t \geq \delta)
   \leq 4 \P (B_{6 \ee} \geq \delta) .$$
As before, the choice of $C_4$ implies that
$\P (R_1) = o (\ee^3 )$
and summing the upper bounds on $\P (R_1)$, $\P (R_2)$ and $\P (R_3)$ 
proves the lemma.   
\endpf

\noindent{\sc Proof of  Proposition}~\ref{pr:2 dim}: The
inequality~(\ref{eq:quadrants}), and the theorem, follow directly
from Lemmas~\ref{thirdlemma} and~\ref{fourthlemma}:
for sufficiently small $\ee > 0$
$$
\P (W_1 \in Q_3 | A_{\ee}) 
= {\P (\{ W_1 \in Q_3 \} \cap A_\ee) \over \P (A_\ee)} 
\leq  {C_3 \ee^3 \log (1/\ee)^3 \over C_1 \ee} 
 =  C \ee^2 \log (1/\ee)^3  .
$$
\endpf

\section{Recovery of the endpoint}
\setcounter{equation}{0}
\label{recov}
In this section,
let $(\Omega , \F , \P)$ be the space of continuous functions 
$\omega : [0,1] \to \Rd$, endowed with the Borel $\sigma$-field $\F$ 
(in the topology of uniform convergence) and Wiener measure $\P$ 
on paths from the origin.  
We write simply $\mu$ instead of $\mu_\Theta$ for the occupation
measure of the spherical projection $(\sph(\omega_t), 0 \le t \le 1)$.
So $\mu$ is a measurable map from $(\Omega , \F)$ to the space
$(\omxxx , \F_2)$ of Borel probability measures on $S^{d-1}$.
For a subinterval $I$ of $[0,1]$, say $I = (a,b)$ or $I = [a,b]$, 
let $\omega I$ denote the range of the restriction of $\omega$ to $I$.

We will use some known topological facts about Brownian motion in
dimensions $d \geq 3$:  
\begin{enumerate}
\item[(1)] If $I$ and $J$ are disjoint open sub-intervals of $[0,1]$, then
$\P$ almost surely the random set 
$ \{ \sph(\omega_t), t \in I \}$ does not contain 
$ \{ \sph(\omega_t), t \in J \}$.
\item[(2)] Almost every Brownian path $\omega : [0,1] 
\rightarrow \Rd$ has a sequence of cut-times $t_n \uparrow 1$, that is,
$\omega (0,t_n) \cap \omega (t_n , 1) = \emptyset$.  
\item[(3)] With probability 1, no cut-point is a double point.  Formally,
for $\P$-almost every $\omega$, if $\omega [0,1] \setminus 
\{ \omega (t) \}$ is not connected, then $\omega (s) \neq \omega (t)$
for $s \neq t$.  
\end{enumerate}
Fact~(1) follows easily from Fubini's theorem.
Fact~(2) is proved in Theorem~2.2
of Burdzy~\cite{Burdzy89} (see also \cite{Burdzy95}) and fact~(3) is
proved in Theorem~1.4 of Burdzy-Lawler~\cite{Burdzy-Lawler90}.
  
The following lemma
contains the probabilistic content of the argument and is proved at
the end of this section.  Facts~(2) and~(3) are true when $d=2$ as
well, which is all that is needed to establish the remark after
Theorem~\ref{th:recover}.

\begin{lmm} {lem:rho}
Let $D$ be a ball in the sphere $S^{d-1} \subseteq \Rd$.  Then there is a 
measurable function $\rho_D : \omxxx \to \R^+$ such that $\P$-almost surely,
\begin{equation} \label{eq:rho}
\rho_D (\mu (\omega)) = \sup \{ |\omega (t)| : \pi (\omega (t)) \in D \} .
\end{equation}
\end{lmm}

To construct $\psi$ from Lemma~\ref{lem:rho} and fact~(1), let $\C_j$
be a finite cover of $S^{d-1}$ by balls of radius $2^{-j}$ and let 
$\cen (D)$ denote the center of the ball $D$.  
The set
of limit points of a sequence $\{ S_n \}$ of elements in $\omxxxx$,
defined by $\{ x : \liminf_n d(x , S_n) = 0 \}$, is called the Hausdorff
limsup, denoted $\limsup_{n \to \infty} S_n$.
Observe that if $S_n$ are $\omxxxx$-valued random variables, then 
$\limsup_{n \to \infty} S_n$ is measurable as well.  

\begin{lmm} {lem:psi}
Define measurable functions $A_j : \omxxx \to \omxxxx$ to be the 
sets of vectors
$$A_j (\mu) := \{ \rho_D (\mu) \cen (D) : D \in \C_j \} .$$
Then $\psi := \limsup_{j \to \infty} A_j$ satisfies~(\ref{eq:psi}):
$$\psi \circ \mu (\omega) = \range (\omega) \hspace{.5in} 
   \mbox{for almost every } \omega .$$
\end{lmm} 

\noindent{\sc Remark:} 
In fact, from the proof we see that $\psi = 
\lim A_j$ almost surely when $\mu = \mu (\omega)$ and $\omega$ is
chosen from $\P$.  

\noindent{\sc Proof:} It is easy to see that $\psi \circ \mu (\omega) 
\subseteq \range (\omega)$ for every $\omega$:
if $D \in \C_j$ then $\rho_D (\mu) (\omega) \cen (D)$ is equal to 
$|\omega (t)| \cen (D)$ for some $t$ with $|\pi (\omega (t)) - \cen (D)| 
\leq 2^{-j}$, and is hence within $2^{-j}|\omega (t)|$
 of the point $\omega (t) \in \range (\omega)$; 
since $\range (\omega)$ is closed and $j$ is 
arbitrary, all limit points of sequences $\{ x_j  \}$ with $x_j \in A_j$ 
are in $\range (\omega)$.  

To see that $\range (\omega) \subseteq \psi \circ \mu (\omega)$, 
fix  $t \in (0,1)$  and consider $x =\omega(t) \in \range(\omega)$.
For any $\eps>0$, choose a 
$\delta > 0$ such that $|\omega (s) - \omega (t)| \le \eps$ when 
$|s - t| \le \delta$, and $|\omega (s)|<|x|$
for $0 \le s \le \delta$.  By fact~(1), 
the union $\pi(\omega[\delta,t-\delta]) \cup \pi (\omega [t + \delta , 1])$ 
does not cover $\pi (\omega [t - \delta , t + \delta])$.  
Thus we may choose an 
open ball $D$ intersecting $\pi (\omega [t - \delta , t + \delta])$ 
such that $\pi (D)$ is disjoint from 
$\pi(\omega[\delta,t-\delta]) \cup \pi (\omega[t+\delta , 1])$.  
For any $D' \subseteq D$, it follows that $|\rho_{D'} (\mu) - |\omega (t)|| 
\leq \eps$.  For sufficiently large $j$ there is a ball 
$D' \in \C_j$ with $x \in D' \subseteq D$, which implies $A_j$ 
contains a point $\rho_{D'} (\mu) \cen (D')$ within $2^{-j}|x|+\eps$ of $x$.  
Since $\eps$ and $j$ are arbitrary, $x$ is a limit point of the sets $A_j$.   
\endpf

The construction of $\varphi$ from here uses two further non-probabilistic
lemmas.

\begin{defn}
Define the map $N_\delta : \omxxxx \to \{ 0 , 1, 2 , \ldots , \infty \}$ 
by setting $N_\delta (S)$ to be the number $N$ of connected components of 
the closed set $S$ that have diameter at least $\delta$. 
\end{defn}

\begin{lmm} {lem:N}
For each $\delta$ the map $N_\delta$ is measurable.  
\end{lmm}
\noindent{\sc Proof}: It suffices to show this when $S$ is a subset of 
the unit ball.  It will be convenient to have a nested sequence of sets 
$\grid_1 \subseteq \grid_2 \subseteq \cdots$ such that $\grid_j$ is 
$2^{-j-1}$-dense in the unit ball.  (To construct this, inductively 
choose $\grid_j$ to be a maximal set with no two points within
distance $2^{-j-1}$.) The sets
$\balls_j$ defined to be the set of balls of radius $2^{-j}$ centered 
at points of $\grid_j$, form a sequence of covers of the unit ball
such that each element of $\balls_{j+1}$ is contained in an element of
$\balls_j$.

For each $j$ and each $S \in \omxxxx$ let 
$$X_j (S) = \bigcup \{ D \in \balls_j : D \cap S \neq \emptyset \}.$$  
Let $P_j$ be the set of 
connected components of $X_j (S)$ viewed as subsets of $\balls_j$.
In other words, $P_j (S) = \{ \C \subseteq \balls_j : \bigcup \C
\mbox{ is a component of } X_j (S) \}$.  By the finiteness of
$\balls_j$, we see that each $P_j$ is measurable.  Since each
$D \subseteq X_j (S)$ is contained in a ball $D' \in \balls_{j-1}$
also intersecting $S$, $X_j \subseteq X_{j-1}$ and hence each component
of $X_j$ is contained in a unique component of $X_{j-1}$.  This
defines a map $\edge_j: P_j \to P_{j-1}$ which is measurable since 
it depends only on $P_j$ and $P_{j-1}$.  Letting $P_{j , \delta}$
be the subset of $P_j$ consisting of components of diameter at least
$\delta$, it is clear that $\edge_j$ maps $P_{j , \delta}$ to
$P_{j-1 , \delta}$ and that these are measurable.  

\noindent{\bf Claim:} $N_\delta (S)$ is the cardinality of the inverse 
limit of the system $\{ P_{j , \delta} , \edge_j : j \geq 1 \}$.  Indeed, 
suppose that $\{ x_j^{(i)} \}$ satisfy $x_j^{(i)} \in P_{j , \delta}$ and 
$\edge_j (x_j^{(i)}) = x_{j-1}^{(i)}$ for all $j$ and $i = 1 , 2$.  Letting
$\Set (x_j^{(i)}) := \bigcup x_j^{(i)}$ denote the set of points in the
component $x_j^{(i)}$, we see that $\bigcap_j (\Set (x_j^{(i)}))$ 
are non-empty subsets of $S$ and lie in different 
components unless $x_j^{(1)} = x^{(2)}_j$ for all $j$.  Conversely,
if $x$ and $y$ are points of $S$ lying in different
connected components, then $S$ is contained in a disjoint union 
$X^\ee \cup Y^\ee$ for some sets $X , Y$ with $x \in X$, $y \in Y$
(where $Z^\ee$ denotes the set of points within $\ee$ of the set $Z$).
It follows that for each $j$ there is an $x_j \in P_{j , \delta}$ with 
$x \in \bigcup x_j$, there is a $y_j \in P_{j , \delta}$ with 
$y \in \bigcup y_j$, 
and that for $2^{-j} < \ee$, $x_j \neq y_j$.  

Finally, the cardinality of the inverse limit is easily seen to be 
measurable.  Say $x_j \in P_{j , \delta}$ is a survivor if for each $k > j$
there is some $y_k \in P_{k ,\delta}$ with $\bigcup y_k \subseteq \bigcup x_j$.
The set of survivors is clearly measurable, and the cardinality
of the inverse limit is the increasing limit of the number of survivors
in the set $P_{j , \delta}$ as $j \to \infty$.   
\endpf

The endpoint $\omega(1)$ will be recovered from $\omega[0,1]$ 
as the only nonzero limit point of cutpoints, which is not a cutpoint
itself. To justify measurability of this operation,
the following definition and lemma are useful.  
\begin{defn} 
Let $\cut_\delta (S)$ denote the set of $\delta$-cutpoints of $S$, 
that is, those $x \in S$ such that $S \setminus x$ has at least 
two components of diameter at least $\delta$ (note: if $S$ is not 
connected this may be all of $S$).  
For each positive integer $j$ and each $\delta > 0$, define 
the measurable function $A_{\delta , j} : \omxxxx \to \omxxxx$ by 
$$A_{\delta , j} (S) := \bigcup \{ D' \in \balls_j : D' \cap S
   \neq \emptyset \mbox{ and } N_\delta (S \setminus D') \geq 2 \} .$$
Let
$$A_\delta := \limsup_{j \to \infty} A_{\delta , j} .$$
\end{defn}

\begin{lmm} {lem:delta-cut}
Let $f : [0,1] \to \Rd$ be any continuous function and denote its range
by $S$.  Fix $\delta' > \delta > \delta'' > 0$.  Then
\begin{equation} \label{eq:delta-cut} 
\cut_{\delta'} (S) \subseteq A_\delta \subseteq \cut_{\delta''} (S) .
\end{equation}
\end{lmm}

\noindent{\sc Proof:} Suppose first that $x$ is a $\delta'$-cutpoint of $S$.  
Let $T$ and $U$ be two components of $S \setminus x$ of diameter at least 
$\delta$.  If $D$ is a ball of radius $\ee < (\delta' - \delta) / 2$ 
containing $x$, then $S \setminus D$ will have at least two components 
of diameter at least $\delta' - 2 \ee$.  Thus $x \in A_{\delta , j}$ for
$2^{-j} < \ee$, hence $x \in A_\delta$.

Suppose now that $x \in A_\delta$ and let $\{ D_n \}$ be balls converging
to $x$ in the Hausdorff metric, such that each intersects $S$ and
has $N_\delta (S \setminus D_n) \geq 2$.  Let $\{ D_n' \} $ be balls 
with diameters going to zero such that $\bigcup_{j=n}^\infty D_j 
\subseteq D_n'$.  Then $N_{\delta''} (S \setminus D_n') \geq 2$ when 
$n$ is large enough so that the diameter of $D_n'$ is at most $\delta 
- \delta''$.  
\begin{quote}
Claim: there are points $x_1 , \ldots , x_k$ and an $N_0$ such that for
$n \geq N_0$, each component of $S \setminus D_n'$ of diameter at least
$\delta''$ contains one of $x_1 , \ldots , x_k$.  

Proof: Pick $N_0$ so that $D_n \subseteq \ball (x , \delta'' / 2)$ 
when $n \geq N_0$.  Pick $\ee > 0$ such that $|f(s) - f(t)| < \delta'' / 2$
when $|s - t| \leq \ee$.  The open set $\{ t : |f(t) - x| > \delta'' / 2 \}$
decomposes into a countable set of intervals.  
At most $k := \lfloor 1 / \ee \rfloor$ of these intervals $(u_j , v_j) ,
j = 1 , \ldots , k$ can have $v - u \geq \ee$, and these are the 
only ones containing times $t$ with $|f(t) - x| \geq \delta''$.  
Since $S$ is connected, every component
$G$ of $S \setminus D_n'$ intersects $\partial D_n'$, and if $G$ has
diameter at least $\delta''$ then $G$ must contain one of the $k$ 
sojourns $f(u_j , v_j)$.  Choose $x_j \in f(u_j , v_j)$.  
\end{quote}

Since $N_{\delta''} (S \setminus D_n') \geq 2$ for all $n \geq N_0$, 
there are $i < j \leq k$ such that infinitely many of the sets
$S \setminus D_n$ have distinct components $G_n$ and $H_n$ of size 
at least $\delta''$ containing $x_i$ and $x_j$ respectively.  
The increasing limits $\bigcup G_n$ and $\bigcup H_n$ must then be
contained in distinct components of $S \setminus \{ x \}$, showing
that $x \in \cut_{\delta''} (S)$.   
\endpf

\noindent{\sc Proof of Theorem~\protect{\ref{th:recover}} 
assuming~Lemma}~\ref{lem:rho}:  Clearly the
sets $A_\delta$ increase as $\delta \to 0$.  Define 
$$\varphi(S) = (\limsup_{\delta \to 0} A_\delta) \setminus (\bigcup_\delta 
   A_\delta \cup \{ 0 \}) .$$
We have shown that $\psi \circ \mu (\omega) = \omega [0,1]$ almost surely with
respect to $\P$, and it follows from Lemma~\ref{lem:delta-cut}
that $\varphi (S) \cup \{ 0 \}$ is the topological boundary of the set 
of cut-points of $S$.  Fact~(2) then implies that $\omega (1) \in \varphi(S)$.
On the other hand, let $x = \omega (t)$ be any limit of cut-points, where 
$0 < t < 1$.  Thus there are times $t_j \to t$ with $\omega (t_j)$
a cut-point. By fact (3), the sets
$\omega (t_j , 1)$ and $\omega (0 , t_j)$ are disjoint,
and each of them is connected.
For $t_j > t$, the set $\omega (t_j , 1)$ is disjoint
from $\omega (0 , t)$ so if $t_j \downarrow t$, then $\omega (t,1)$
is disjoint from $\omega (0 , t)$.  Likewise if $t_j \uparrow t$ then
$\omega (0,t)$ is disjoint from $\omega (t , 1)$, hence $t$ is a 
cut-time.  This shows that $x \notin \varphi(S)$, so the only limits of 
cut-points that are not cutpoints are $\omega (0)$ and $\omega (1)$,
which completes the proof.   
\endpf

To prove Lemma~\ref{lem:rho} we state several more lemmas.  The cases
$d = 3$ and $d \geq 4$ differ slightly in that the estimates required
for two-dimensional balls ($d = 3$) include logarithmic terms.  
Since recovery of the endpoint in dimension $d \ge 4$ can be
reduced to the three-dimensional case, and since the estimates for
two-dimensional balls are strictly harder than for higher-dimensional
balls, we assume for the remainder of the proof that $d = 3$. 
The formula for $\rho_D$ in this case is given by: 
\begin{equation} \label{eq:def-rho}
\rho_D (\mu) := \left [\limsup_{D' \subseteq D , 
   \rad(D') \to 0} {\mu (D') \over 2 \rad (D')^2 
   \log^2 \rad (D')} \right ]^{1/2} .
\end{equation}
We remark that when $d > 3$, the term $\log^2 \rad (D')$ is replaced
by $|\log \rad (D')|$ and the constant 2 in the denominator changes 
as well; this is due to the different normalization needed for ``thick points''
in dimension $3$ and higher, see \cite{DPRZ99a}.

We begin by quoting two results from 
Dembo, Peres, Rosen and Zeitouni \cite{DPRZ99b}.  

\begin{lmm} {lem:DPRZ1}(\cite{DPRZ99b}, Theorem 1.2).
Let $(W_t : t \geq 0)$ be a standard Brownian motion in $\R^2$.
Let $\rad (D)$ denote the radius of the ball $D$.  Then for any
fixed $A > 0$, 
\begin{equation} \label{eq:d=2}
\limsup_{\rad (D) \to 0} {\int_0^A \one_D (W_t) dt \over 
   \rad (D)^2 \log^2 \rad (D)} = 2 \quad a.s.
\end{equation}
\end{lmm}
\endpf

\begin{lmm} {lem:DPRZ2} (\cite{DPRZ99b}, Lemma 2.1). 
Let $Z_t = \int_0^t \one_D (W_t) \, dt$ be the occupation time of a 
standard two-dimensional Brownian motion up to time $t$ in a ball
$D$ of radius $r$.  Then for each $t > 0$ there is some $\lambda > 0$
not depending on $r$ for which $\E e^{\lambda Z_t / (r^2 |\log r|)} 
< \infty$.  Consequently, $\P (Z_t > A r^2 \log (1/r) < C e^{-\gamma A}$
for some positive $C$ and $\gamma$.
\end{lmm}

\noindent{\sc Proof:} Dembo et al prove the result when the Brownian
motion is started at radius $r$ (in their notation $r = r_1 = r_2$)
and the time $t$ is instead the time to hit a ball of fixed radius
$r_3 = O(1)$.  Accomodating these changes is trivial.  
\endpf

We now state three more lemmas which together imply Lemma~\ref{lem:rho}.

\begin{lmm} {lem:eps-1}
Let $(W_t : t \geq 0)$ be a standard three-dimensional Brownian motion.
For $0 \leq a < b \leq 1$, let $\mu_{a,b}$ be projected occupation 
measure in the time interval $[a,b]$, i.e., for $D \subseteq S^2$,
$$\mu_{a,b} (D) := \int_a^b \one_D (\pi (W_t)) \, dt .$$
Then for each ball $D \subseteq S^2$ and each $\ee > 0$, with probability 1,
\begin{equation} \label{eq:eps-1}
\limsup_{D' \subseteq D , \rad(D') \to 0} {\mu_{\ee , 1} (D') 
   \over \rad (D')^2 \log^2 \rad (D')} \leq 
   2 (\sup \{ |W_t| : \pi (W_t) \in D \})^2 .
\end{equation}
\end{lmm}

\begin{lmm} {lem:0-eps}
In the notation of the previous lemma, there is a constant $c_2$ 
such that for each $\ee > 0$, with probability 1,
\begin{equation} \label{eq:0-eps}
\limsup_{D' \subseteq D , \rad(D') \to 0} {\mu_{0 , \ee} (D') 
   \over \rad (D')^2 \log^2 \rad (D')} \leq 
   c_2 (\sup \{ |W_t| : t \in [0,\ee] \})^2 .
\end{equation}
\end{lmm}

\begin{lmm} {lem:lower}
For each $t \in (0,1)$, with probability 1,
$$\limsup_{D \to \pi (W_t)} {\mu (D) \over \rad (D)^2 
   \log^2 \rad (D)} \geq 2 |W_t|^2 . $$
\end{lmm}

To see why Lemma~\ref{lem:rho} follows from Lemmas~\ref{lem:eps-1}~-
\ref{lem:lower}, define $\rho_D$ as in equation~(\ref{eq:def-rho}).
Since the limsup may be taken over balls with rational centers and radii,
$\rho_D$ is measurable.  Lemmas~\ref{lem:eps-1} and~\ref{lem:0-eps}
together imply that with probability 1, for all $\ee > 0$,
$$\rho_D (\mu (\omega)) \leq \left [ (\sup \{ |W_t| : \pi (W_t) \in D 
   \})^2 + (c_2 / 2) (\sup \{ |W_t| : t \in [0,\ee] \})^2
   \right ]^{1/2} , $$
and sending $\ee$ to 0 shows that the LHS of~(\ref{eq:rho}) is less than or
equal to the RHS.  On the other hand, applying Lemma~\ref{lem:lower}
for all rational $t$ shows that with probability 1,
$$\rho_D (\mu (\omega)) \geq \sup \{ |\omega (t)| : \pi (\omega (t)) \in 
   \INT (D) , t \mbox{ rational} \} $$
which yields the reverse inequality.  It remains to 
prove Lemmas~\ref{lem:eps-1}~-~\ref{lem:lower}.

\noindent{\sc Proof of Lemma}~\ref{lem:eps-1}:
Covering $D$ with small balls, it suffices to assume $\rad (D) < \delta$
and prove an upper bound of $(1 + o(1))$ times the RHS of~(\ref{eq:eps-1})
as $\delta \to 0$.  Let $\beta : S^2 \to \R^2$ be a conformal map
with Jacobian going to 1 near $\cen (D)$.  For example, take $\beta$
to be stereographic projection from the antipode to $\cen (D)$ to a plane
(identified with $\R^2$) tangent to $S^2$ at $\cen (D)$.  
The path $\{ \pi (W_t) : t \geq \ee \}$ is a time-changed 
Brownian motion on $S^2$, and in particular, $\pi (W_{G(t)}))$ is 
a Brownian motion started from $\pi (W_1)$, where $G(t)$ is defined
by $\int_{G(t)}^1 |W_s|^{-2} \, ds = t$.
Similarly, $( X_t := \beta (\pi (W_{G(H(t))})) , t \in 
[0 , M := H^{-1} (G^{-1} (\ee))] )$ is a Brownian motion in $\R^2$,
where $M$ is random and $H(t)$ is another time change, with $|H'|$ 
going to 1 uniformly as $\rad (D) \to 0$ and $\pi (W_{G(H(t))})$ is in $D$.  

Let $D'$ be any ball inside $D$.  Let $D''$ be a ball containing 
$\beta (D')$ and observe that we can take $\rad (D'') / \rad (D') \to 1$
uniformly over $D' \subseteq D$ as $\rad (D) \to 0$.  When $\pi (
W_s) \in D$, $G' (s) \geq \sup \{ |W_t| : \pi (W_t) \in D \}^2$.  Thus 
\begin{eqnarray*}
\mu_{\ee , 1} (D') & = & \left | \left \{ t \in [\ee , 1] : \pi (W_t)
   \in D' \right \} \right | \\[2ex]
& \leq & \left | \left \{ G(H(s)) : \beta (\pi (W_{G(H(s))})) \in D''
   \right \} \right | \\[2ex]
& \leq & \sup \{ |W_t| : \pi (W_t) \in D \}^2 \sup H' \left | \left \{
   s :  \beta (\pi (W_{G(H(s))})) \in D'' \right \} \right | \\[2ex]
& = & \sup \{ |W_t| : \pi (W_t) \in D \}^2 \sup H' \left | \left \{
   s :  X_s \in D'' \right \} \right | \\[2ex]
& \leq & (2 + o(1)) \sup \{ |W_t| : \pi (W_t) \in D \}^2 
   r(D'')^2 \log^2 ({1 \over r(D'')}) 
\end{eqnarray*}
by Lemma~\ref{lem:DPRZ1} and the convergence of $H'$ to 1.  
\endpf

\noindent{\sc Proof of Lemma}~\ref{lem:0-eps}: Let $D$ be any ball in
$S^2$ with center $x$.  Let $\beta_x$ be projection to the orthogonal
complement of $x$ in $\R^3$.  If $\pi (W_t) \in D$ then $\beta_x (W_t) 
\in \ball (0 , s)$ for $s := \rad (D) \sup \{ |W_t| : t \in [0,1] \}$.  
For fixed $x$, $\beta_x (W_t)$ is a standard Brownian motion, so
an application of the Lemma~(\ref{lem:DPRZ2}) yields 
\begin{eqnarray*}
&& \P \left ( {\mu_{0,\ee} (D) \over \rad (D)^2 \log^2 \rad (D)} 
   \geq c (\sup \{ |W_t| : t \in [0,1] \} )^2 \right ) \\[2ex]
& \leq & \P \left ( {\int_0^1 dt \, \one_{\ball (0 , s)} 
   (\beta_x (W_t)) \over s^2 |\log s| |\log r(D)| (\log (\rad (D)) / \log s)} 
   \geq c \right ) \\[2ex]
& \leq & C \rad (D)^{- \gamma c \log \rad (D) / \log s}  .
\end{eqnarray*}
We may choose $c_2$ so that $c_2 \gamma > 2$, and find classes $\C_r$ 
of balls of radius $r$ so that for any $\ee > 0$, for sufficiently 
small $r$, any ball of radius $(1 - \ee) r$ is contained in some 
element of $\C_r$.  One can arrange for $|\C_r| = O(1/r)^{c_2 c_0 - \delta}$,
where $c_2 c_0 - \delta > 2$, ensuring that 
$$\P \left (\exists D \in \C_r : {\mu_{0,\ee} (D) \over \rad (D)^2 
   \log^2 \rad (D)} \geq c_2 (\sup \{ |W_t| : t \in [0,1] \} )^2 \right ) 
   = o(r^\delta) .$$
Summing over $r = (1 - \alpha)^n$ and using Borel-Cantelli shows that 
the limsup on the LHS of~(\ref{eq:0-eps}) is at most $(1 - \alpha)^{-2}
c_2 \sup \{ |W_t| : t \in [0,\ee] \}^2$, proving the lemma since $\alpha$
may be chosen arbitrarily small.
\endpf

\noindent{\sc Proof of Lemma}~\ref{lem:lower}: Fix $t \in (0,1)$.
Define $\beta , G$ and $H$ as in the proof of Lemma~\ref{lem:eps-1},
so that $( X_s := \beta (\pi (W_{G(H(s))})) )$ is a planar Brownian
motion.  For any $\ee > 0$, Lemma~\ref{lem:DPRZ1} yields a random sequence 
of balls $D_n \to 0$ in $\R^2$ with
$$\int_0^M {\one_{D_n} (X_s) \, ds \over \rad (D_n)^2 \log^2 \rad (D_n)}
   \to 2 .$$
With probability 1, $W_t$ is a single value, i.e., $W_t \neq W_s$
for $t \neq s$, in which case for $n$ sufficiently large, $X_s \in D_n$
implies $|W_{G(H(s))}| \to |W_t|$ and $G(H(s)) \to t$.  The sets 
$\beta^{-1} (D_n)$ are contained in balls $D_n'$ with $\rad (D_n')
/ \rad (D_n) \to 1$, so 
$$\int_0^M {\one_{D_n'} (\pi (W_{G(H(s))}) \, ds \over \rad (D_n')^2 
   \log^2 \rad (D_n')} \to 2 .$$
Changing variables reduces this integral to
$$\int_\ee^1 {\one_{D_n'} (\pi (W_u)) (G \circ H)' ((G \circ H)^{-1} (u))
   \, du \over \rad (D_n')^2 \log^2 \rad (D_n')} $$ 
and since $(G \circ H)' = (1 + o(1)) |W_t|^{-2}$ uniformly on an interval
containing $H^{-1} (G^{-1} (t))$,  we get 
$$|W_t|^{-2} {\mu_{\ee , 1} (D_n') \over \rad (D_n') \log^2 \rad (D_n')}
   \to 2 ,$$
proving the lemma.  
\endpf

\noindent{\bf Acknowledgments:} We are grateful
to Jason Schweinsberg for suggesting the example in 
 Proposition~\ref{pr:2 dim}. We thank
G\'erard Letac and Zhan Shi for the reference \cite{purirubin70} to
the characterization of the exponential distribution used in the proof of
Proposition \ref{pr:characterize}.  We thank Tom Salisbury and the referee
for useful comments.  We are indebted to MSRI and 
to the organizers of the 1998 program on stochastic analysis there, for
the chance to join forces.

\end{document}